%% file: softSquare14A.tex
\documentclass[12pt,oneside,english]{amsart}
\usepackage{ae,aecompl}

\usepackage[T1]{fontenc}
\usepackage[latin9]{inputenc}
\usepackage[letterpaper]{geometry}
\usepackage{array}
\usepackage{amsthm}
\usepackage{amssymb}

\makeatletter

\providecommand{\tabularnewline}{\\}

\numberwithin{equation}{section} %% Comment out for sequentially-numbered
\numberwithin{figure}{section} %% Comment out for sequentially-numbered
\theoremstyle{plain}
\newtheorem{thm}{Theorem}[section]
  \theoremstyle{plain}
  \newtheorem{lem}[thm]{Lemma}

\usepackage{hyperref}

\usepackage[all]{xy}

\usepackage{mathbbol}

\usepackage{ifthen}

\newboolean{Details}

\setboolean{Details}{false}
\usepackage{changebar}

\makeatother

\usepackage{babel}

\begin{document}

\title{Noncommutative Semialgebraic sets and Associated Lifting Problems}

\author{Terry A. Loring}

\author{Tatiana Shulman}

\address{Department of Mathematics and Statistics, University of New Mexico,
Albuquerque, NM 87131, USA.}

\address{Department of Mathematics, University of Copenhagen, Universitetsparken
5, DK-2100 Copenhagen \O , Denmark}



\maketitle
\markleft{Noncommutative Semialgebraic sets and Associated Lifting Problems (version 13)}\markright{Noncommutative Semialgebraic sets and Associated Lifting Problems (version 13)}

\begin{center}
\ifthenelse{\boolean{Details}}{\cbstart  {\Large Fully Detailed Version}\cbend}{}
\par\end{center}
\begin{abstract}
We solve a class of lifting problems involving approximate polynomial
relations (soft polynomial relations). Various associated $C^{*}$-algebras
are therefore projective. The technical lemma we need is a new manifestation
of Akemann and Pedersen's discovery of the norm adjusting power of
quasi-central approximate units. 

A projective $C^{*}$-algebra is the analog of an absolute retract.
Thus we can say that various noncommutative semialgebraic sets turn
out to be absolute retracts. In particular we show a noncommutative
absolute retract results from the intersection of the approximate
locus of a homogeneous polynomial with the noncommutative unit ball.
By unit ball we are referring the $C^{*}$-algebra of the universal
row contraction. We show projectivity of alternative noncommutative
unit balls. 

Sufficiently many $C^{*}$-algebras are now known to be projective
that we are able to show that the cone over any separable $C^{*}$-algebra
is the inductive limit of $C^{*}$-algebras that are projective.
\end{abstract}
\ifthenelse{\boolean{Details}}{\cbstart  Bars on the right indicate
details to be omitted when we circulate this. \cbend}{}

\tableofcontents{}

\section{Introduction\label{sec:Introduction}}

Lifting problems for relations in $C^{*}$-algebras have tended to
have \emph{ad hoc} solutions. Olsen and Pedersen prove in \cite{OlsenPedLifting}
that a nilpotent aways has a nilpotent lift, specifically that given
$x$ in a $C^{*}$-algebra quotient $A/I$ with \[
x^{n}=0\]
there is always $X$ in $A$ with $\pi(X)=x$ and $X^{n}=0.$ Their
proof is rather different from the techniques Akemann and Pedersen
used in \cite{AkePederIdealPerturb} to show that for $x$ in $A/I$
with \[
\left\Vert x^{n}\right\Vert \leq\epsilon\quad(\epsilon>0)\]
there is always a lift $X$ with $\left\Vert X^{n}\right\Vert \leq\epsilon.$
Different still are the techniques used in \cite{LorPederProjTrans}
to show that the relations describing $\mathbf{C}\mathbf{M}_{n}=C_{0}\left((0,1],\mathbf{M}_{n}\right)$
are liftable: given \textbf{$x_{1},\ldots,x_{n}$ }in $A/I$ satisfying
all the relations\begin{align*}
\left\Vert x_{j}\right\Vert  & \leq1\quad(\forall j)\\
x_{j}^{*}x_{k} & =0\quad(j\neq k)\\
x_{j}^{*}x_{j} & =x_{k}^{*}x_{k}\quad(\forall j,\forall k)\\
x_{1}x_{1}^{*} & =x_{1}^{*}x_{1}\end{align*}
there are lifts $X_{1},\ldots,X_{n}$ in $A$ that also satisfy these
relations. More recently, $M$-ideals showed up in \cite{Shulman_nilpotents}
to settle the lifting problem for the relations\begin{align*}
\left\Vert x\right\Vert  & \leq1\\
x^{n} & =0.\end{align*}

The lifting results above (most of them, anyway) show various $C^{*}$-algebras
are projective. Projectivity was introduced by Effros and Kaminker,
in \cite{EffrosKaminkerShape}. A $C^{*}$-algebra $P$ is \emph{projective}
if the map\[
\rho\circ\mbox{--}:\hom(P,B)\rightarrow\hom(P,C)\]
is onto whenever $\rho:B\rightarrow C$ is onto.

Projectivity was shown by Blackadar in \cite{Blackadar-shape-theory}
to be the noncommutative analog of a space being an absolute retract
(AR). The analog of absolute neighborhood retract is semiprojectivity,
which we will not discuss in detail in this paper except in Section~\ref{sec:Soft-Cylinders}.

Systematic investigations of projectivity exist, but only in the case
of at-most one-dimensional spectrum. There was a study of $C_{0}(X)$
for $X^{+}$ a tree in \cite{LoringProjectiveCstar}. Chigogidze and
Dranishnikov solved the general question for $C_{0}(X)$ being projective,
in \cite{chigogidzeNonComARs}. The answer is that $C_{0}(X)$ is
projective if and only if $X^{+}$ is a dedrite. $ $The finite mapping
telescopes associated to inclusions of finite-dimensional $C^{*}$-algebras
were shown to be projective, in \cite{LorPederProjTrans}. In the
later terminology of \cite{ELP-anticommutation}, this says we have
projectivity for a large class of one-dimensional noncommutative CW
complexes.

``NC'' will stand for ``noncommutative.'' Thus noncommutative CW complex
becomes NCCW.

In the commutative case, very sweeping statements can be made about
what spaces are AR or ANR. For example, every compact semialgebraic
set in finite-dimensional Euclidean space is an absolute neighborhood
retract. See \cite[p.\ 79]{VassilievAppliedPicardLef} and \cite{LojasiewiczSemiAn}
for precise results and definitions. A subset of Euclidean space is
said to be \emph{semialgebraic} if it is the union of solution sets
of polynomial equations and polynomial inequalities. As we are interested
in closed and connected sets, it will suffice to have in mind sets
of the form\[
\left\{ \left(x_{1},\ldots,x_{n}\right)\in\mathbb{R}^{n}\left|\strut\, p_{j}\left(x_{1},\ldots,x_{n}\right)\leq\epsilon_{j}\mbox{ for }j=1\ldots J\right.\right\} ,\]
where the $p_{j}$ are polynomials.

This general result about semialgebraic sets being ANR cannot translate
directly to $C^{*}$-algebras. We know that for the unit disk $\mathbb{D},$
the $C^{*}$-algebra \[
C_{0}\left(\mathbb{D}\setminus\{0\}\right)\cong C^{*}\left\langle x\left|\,\strut x^{*}x=xx^{*},\ \left\Vert x\right\Vert \leq1\right.\right\rangle \]
fails to be projective. Normals don't generally lift to normals. Some
normals fail to have partial lifts, and are bounded away from other
normals that have partial lifts. To get technical, $C_{0}\left(\mathbb{D}\setminus\{0\}\right)$
is not even weakly semiprojective (\cite{EilersLoringContingenciesStableRelations}).

There is a way to avoid the difficulty posed by this nonliftable example
other than keeping to small dimension. We will \emph{avoid exact relations}. 

\begin{table}
\begin{tabular}{|>{\centering}m{0.1\textwidth}|>{\centering}m{0.16\textwidth}|>{\centering}m{0.25\textwidth}|>{\centering}m{0.16\textwidth}|>{\centering}m{0.18\textwidth}|}
\hline 
\multicolumn{1}{|>{\centering}p{0.11\textwidth}|}{Gen\-erators} & Individual restrictions & Other Relations & Name or Comment & Credit\tabularnewline
\hline
\hline 
$x$ & \vspace*{0.1cm}
\begin{minipage}[t]{1\columnwidth}%
$\left\Vert x\right\Vert \leq1$

$x^{n}=0$%
\end{minipage}\vspace*{0.1cm}
 &  &  & Shulman \cite{Shulman_nilpotents}(Olsen, Pedersen, Loring for $\left\Vert x\right\Vert <1.$)\tabularnewline
\hline 
$x$ & \vspace*{0.1cm}
\begin{minipage}[t]{1\columnwidth}%
$\left\Vert x\right\Vert \leq1$

$\left\Vert x^{n}\right\Vert \leq C$%
\end{minipage}\vspace*{0.1cm}
 &  &  & Akemann, Pedersen \cite{AkemannLeftIdeal}\tabularnewline
\hline 
$x,y$ & %
\begin{minipage}[t]{1\columnwidth}%
$\left\Vert x\right\Vert \leq1$

$\left\Vert y\right\Vert \leq1$%
\end{minipage} & \vspace*{0.1cm}
\begin{minipage}[t]{1\columnwidth}%
$x^{*}x=y^{*}y$

$x^{*}y=y^{*}x=0$

$x^{2}=y^{2}=0$%
\end{minipage}\vspace*{0.1cm}
 & $C_{0}\left((0,1],\mathbf{M}_{3}\right)$ & Loring, Pedersen \cite{LorPederProjTrans}\tabularnewline
\hline 
$x,y$ & \vspace*{0.1cm}
\begin{minipage}[t]{1\columnwidth}%
$-1\leq x\leq1$ 

$-1\leq y\leq1$%
\end{minipage}\vspace*{0.1cm}
 & %
\begin{minipage}[t]{1\columnwidth}%
$xy=0$%
\end{minipage} & $C_{0}\left((0,1],\mathbb{C}^{4}\right)$ & Loring \cite{LoringProjectiveCstar}\tabularnewline
\hline 
$x_{1},\ldots,x_{n}$ & %
\begin{minipage}[t]{1\columnwidth}%
$\left\Vert x_{j}\right\Vert \leq1,\,(\forall j)$%
\end{minipage} & \vspace*{0.1cm}
\begin{minipage}[t]{1\columnwidth}%
$x_{i}^{*}x_{i}=x_{j}^{*}x_{j}\quad(\forall i,\forall j)$

$x_{i}^{*}x_{i}=0\quad(\mbox{if }i\neq j)$%
\end{minipage}\vspace*{0.1cm}
 &  & Loring, Pedersen \cite[Example 3.11]{LorPederProjTrans}\tabularnewline
\hline 
$h,k,x$ &  & \vspace*{0.1cm}
\begin{minipage}[t]{1\columnwidth}%
$hk=0$

$0\leq\left[\begin{array}{cc}
\mathbb{1}-h & x^{*}\\
x & k\end{array}\right]\leq1$%
\end{minipage}\vspace*{0.1cm}
 & %
\begin{minipage}[t]{1\columnwidth}%
Usefully in the\\
 $q\mathbb{C}$ picture of\\
 $K$-theory.%
\end{minipage} & Loring \cite{LoringProjectiveKtheory}\tabularnewline
\hline 
$x,y,z,w$ & %
\begin{minipage}[t]{1\columnwidth}%
$0\leq x\leq1$

$0\leq y\leq1$

$0\leq z\leq1$

$0\leq w\leq1$%
\end{minipage} & \vspace*{0.1cm}
\begin{minipage}[t]{1\columnwidth}%
$xy=0$

$zw=0$

$(1-x)z(1-x)=0$

$(1-x)w(1-x)=0$%
\end{minipage}\vspace*{0.1cm}
 & %
\begin{minipage}[t]{1\columnwidth}%
$C_{0}(X)$ where\\
 $X^{+}$ is a tree\\
with four\\
 edges.%
\end{minipage} & Loring \cite{LoringProjectiveCstar}\tabularnewline
\hline 
$x_{1},\ldots,x_{r}$ &  & \vspace*{0.1cm}
\begin{minipage}[t]{1\columnwidth}%
${\displaystyle \left\Vert \sum_{k=1}^{r}x_{k}x_{k}^{*}\right\Vert \leq1}$%
\end{minipage}\vspace*{0.1cm}
 &  & Folklore, functional calculus\tabularnewline
\hline
\end{tabular}

\caption{Some Known Projective $C^{*}$-algebras / Liftable relations\label{tab:Known-Projectives}}

\end{table}

An important instance of Theorem~\ref{thm:MixedVariables} is the
fact that for any \emph{positive }$\epsilon,$ an element $x$ in
a $C^{*}$-algebra quotient $A/I$ with\begin{align*}
\left\Vert x\right\Vert  & \leq1\\
\left\Vert x^{*}x-xx^{*}\right\Vert  & \leq\epsilon\end{align*}
has a lift to $X,$ so $\pi(X)=x,$ with $\left\Vert X\right\Vert \leq1$
and $\left\Vert X^{*}X-XX^{*}\right\Vert \leq\epsilon.$ Put another
way, we show\[
A_{\epsilon}=C^{*}\left\langle x\left|\,\strut\left\Vert x^{*}x-xx^{*}\right\Vert \leq\epsilon,\ \left\Vert x\right\Vert \leq1\right.\right\rangle \]
is projective for all positive $\epsilon.$ Since \[
C_{0}(\mathbb{D}\setminus\{0\})\cong\lim_{\rightarrow}A_{\frac{1}{k}}\]
we have shown $C_{0}(\mathbb{D}\setminus\{0\})$ has a shape system
(c.f.~\cite{Blackadar-shape-theory,BorsukTheoryOfShape,EffrosKaminkerShape})
that is trivial in the sense that all the $C^{*}$-algebras in the
system are projective. It was previously unknown whether $C_{0}(\mathbb{D}\setminus\{0\})$
could be written as an inductive limit of semiprojective $C^{*}$-algebras.

\begin{table}
\begin{tabular}{|>{\centering}m{0.11\textwidth}|>{\centering}m{0.16\textwidth}|>{\centering}m{0.25\textwidth}|>{\centering}m{0.16\textwidth}|>{\centering}m{0.18\textwidth}|}
\hline 
Gen\-erators & Individual restrictions & Other Relations & Name / Remark & Reference\tabularnewline
\hline
\hline 
$h,k$ & \vspace*{0.1cm}
\begin{minipage}[t]{1\columnwidth}%
$-1\leq h\leq1$

$-1\leq k\leq1$%
\end{minipage}\vspace*{0.1cm}
 & %
\begin{minipage}[t]{1\columnwidth}%
$\left\Vert \left[h,k\right]\right\Vert \leq C$%
\end{minipage} & Soft Square & Theorem~\ref{thm:MainTheorem}\tabularnewline
\hline 
$h,k$ & %
\begin{minipage}[t]{1\columnwidth}%
$-1\leq h\leq1$

$-1\leq k\leq1$%
\end{minipage} & \vspace*{0.1cm}
\begin{minipage}[t]{1\columnwidth}%
$\left\Vert \left[h,k\right]\right\Vert \leq C$

$\left\Vert h+ik\right\Vert \leq1$%
\end{minipage}\vspace*{0.1cm}
 & Soft Disk I & Theorem~\ref{thm:MainTheorem}\tabularnewline
\hline 
$x$ & %
\begin{minipage}[t]{1\columnwidth}%
$\left\Vert x\right\Vert \leq1$%
\end{minipage} & \vspace*{0.1cm}
\begin{minipage}[t]{1\columnwidth}%
$\left\Vert \left[x^{*},x\right]\right\Vert \leq2C$%
\end{minipage}\vspace*{0.1cm}
 & Soft Disk I

a second

presentation & \tabularnewline
\hline 
$h,k$ & %
\begin{minipage}[t]{1\columnwidth}%
$-1\leq h\leq1$

$-1\leq k\leq1$%
\end{minipage} & \vspace*{0.1cm}
\begin{minipage}[t]{1\columnwidth}%
$\left\Vert \left[h,k\right]\right\Vert \leq C$

$\left\Vert h^{2}+k^{2}\right\Vert \leq1$%
\end{minipage}\vspace*{0.1cm}
 & Soft Disk II & Theorem~\ref{thm:MainTheorem}\tabularnewline
\hline 
$a,k$ & \vspace*{0.1cm}
\begin{minipage}[t]{1\columnwidth}%
$\left\Vert a\right\Vert \leq1$

$0\leq k\leq1$%
\end{minipage}\vspace*{0.1cm}
 & %
\begin{minipage}[t]{1\columnwidth}%
$\left\Vert \left[a,k\right]\right\Vert \leq C$%
\end{minipage} & Might be useful investigating commutators and square roots & Theorem~\ref{thm:MainTheorem}\tabularnewline
\hline 
$h,p$ & \vspace*{0.1cm}
\begin{minipage}[t]{1\columnwidth}%
$0\leq h\leq1$

$0\leq p\leq1$%
\end{minipage}\vspace*{0.1cm}
 & %
\begin{minipage}[t]{1\columnwidth}%
$\left\Vert h(p^{2}-p)\right\Vert \leq C$%
\end{minipage} &  & Theorem~\ref{thm:MainTheorem}\tabularnewline
\hline 
$h,k,x$ &  & \vspace*{0.1cm}
\begin{minipage}[t]{1\columnwidth}%
$\left\Vert hk\right\Vert \leq C$

$0\leq\left[\begin{array}{cc}
\mathbb{1}-h & x^{*}\\
x & k\end{array}\right]\leq1$%
\end{minipage}\vspace*{0.1cm}
 &  & Theorem~\ref{thm:MainTheorem}\tabularnewline
\hline 
$x,y,z,w$ & %
\begin{minipage}[t]{1\columnwidth}%
$0\leq x\leq1$

$0\leq y\leq1$

$0\leq z\leq1$

$0\leq w\leq1$%
\end{minipage} & \vspace*{0.1cm}
\begin{minipage}[t]{1\columnwidth}%
$\left\Vert xy\right\Vert \leq C$

$\left\Vert zw\right\Vert \leq C$

$\left\Vert (1-x)z(1-x)\right\Vert \leq C$

$\left\Vert (1-x)w(1-x)\right\Vert \leq C$%
\end{minipage}\vspace*{0.1cm}
 &  & Theorem~\ref{thm:SoftTrees}\tabularnewline
\hline 
$x_{1},\ldots,x_{r}$ &  & \vspace*{0.1cm}
\begin{minipage}[t]{1\columnwidth}%
${\displaystyle \left\Vert \sum_{k=1}^{r}\left(x_{k}x_{k}^{*}\right)^{\frac{q}{2}}\right\Vert \leq1}$%
\end{minipage}\vspace*{0.1cm}
 & $1<q<\infty$ & Theorem~\ref{thm:AltNCballs}\tabularnewline
\hline
\end{tabular}

\caption{Some New Projective $C^{*}$-algebras / Lift\label{tab:New-Projectives}able
relations}

\end{table}

An important special case that we study is the approximate zero locus
of a homogeneous NC $*$-polynomial intersected with the NC unit ball.
The homogeneity is imposed to give contractability, and so gives us
an expectation of finding not only semiprojectivity, but projectivity.
By approximate zero locus we mean the universal $C^{*}$-algebra\emph{\[
A_{\epsilon}=C^{*}\left\langle x_{1},\ldots,x_{r}\left|\,\strut\left\Vert p\left(x_{1},\ldots,x_{r}\right)\right\Vert \leq\epsilon,\ \left\Vert \sum_{j=1}^{r}x_{j}x_{j}^{*}\right\Vert \leq1\right.\right\rangle ,\]
}and the ``row contraction'' condition $\sum x_{j}x_{j}^{*}\leq1$
(c.f.~\cite{ArvesonSubalgIII,DavidsonPitts}) is implementing the
intersection with the unit ball. A special case of Theorem~\ref{thm:MixedVariables}
states that $A_{\epsilon}$ is projective for all $\epsilon>0.$

In many cases the relations we can handle have as their universal
$C^{*}$-algebra something that is very unfamiliar. In these cases,
it is perhaps best to see our results as lifting results for the relations
rather than projectivity results for the $C^{*}$-algebras.

Table \ref{tab:Known-Projectives}, lists some known liftable relations.
It is somewhat representative. There are closure results, such as
when $A$ is projective also $\mathbf{M}_{n}(A)$ is projective, which
lead to many more liftable sets of relations, so no table can be complete.

Table \ref{tab:New-Projectives}, lists some of the new examples.
It is not certain these are new projective $C^{*}$-algebras, as projective
$C^{*}$-algebras are contractible and so invariants such as $K$-theory
are of no avail. We can say with some certainty the relations were
not known to lift.

An useful result in topological shape theory is that every compact
metric space is the projective limit of ANRs \cite[IX.1.4]{BorsukTheoryOfShape}.
Blackadar writes in \cite[4.4]{Blackadar-shape-theory}:
\begin{quote}
It is not clear that every C{*}-algebra has a strong shape system
{[}is an inductive limit of semiprojective C{*}-algebras{]}.
\end{quote}
We show in Section~\ref{sec:Cones-are-Limits} that every cone over
a separable $C^{*}$-algebra is the inductive limit of projective
$C^{*}$-algebras.

\section{Quasi-Central Approximate Units Fix Norms\label{sec:Quasi-Central-Approximate-Units}}

Our key tool for lifting is Theorem~\ref{thm:normTechnicalLemma}.
It was extracted from the difficult terrain that is page 127 of Akemann
and Pedersen's paper \cite{AkePederIdealPerturb}.

Approximate units are assumed to satisfy $0\leq u_{\lambda}\le1.$
If $I$ is an ideal in $A$ we let $\pi:A\rightarrow A/I$ denote
the quotient map.
\begin{lem}
\label{lem:approxUnitsLemma} Suppose $I\vartriangleleft A.$ For
any approximate unit $u_{\lambda}$ of $I,$ any $h$ in $A_{+},$
and any real $0\leq\delta\leq1,$\[
\limsup_{\lambda}\left\Vert \left(1-u_{\lambda}\right)^{\frac{1}{2}}h\left(1-u_{\lambda}\right)^{\frac{1}{2}}+(1-\delta)u_{\lambda}^{\frac{1}{2}}hu_{\lambda}^{\frac{1}{2}}\right\Vert \leq\max\left(\left\Vert \pi(h)\right\Vert ,(1-\delta)\left\Vert h\right\Vert \right).\]
\end{lem}
\begin{proof}
We can lift $\pi(h)$ to $k$ with $0\leq k\leq\left\Vert \pi(h)\right\Vert .$
Setting $x=h-k$ we have $x$ in $I$ and\begin{align*}
 & \limsup_{\lambda}\left\Vert \left(1-u_{\lambda}\right)^{\frac{1}{2}}h\left(1-u_{\lambda}\right)^{\frac{1}{2}}+(1-\delta)u_{\lambda}^{\frac{1}{2}}hu_{\lambda}^{\frac{1}{2}}\right\Vert \\
 & \hspace{1in}=\limsup_{\lambda}\left\Vert \left(1-u_{\lambda}\right)^{\frac{1}{2}}k\left(1-u_{\lambda}\right)^{\frac{1}{2}}+(1-\delta)u_{\lambda}^{\frac{1}{2}}hu_{\lambda}^{\frac{1}{2}}\right\Vert .\end{align*}
Now we use the order structure in $A$ and find\begin{align*}
\left(1-u_{\lambda}\right)^{\frac{1}{2}}k\left(1-u_{\lambda}\right)^{\frac{1}{2}}+(1-\delta)u_{\lambda}^{\frac{1}{2}}hu_{\lambda}^{\frac{1}{2}} & \leq\left\Vert k\right\Vert \left(1-u_{\lambda}\right)+(1-\delta)\left\Vert h\right\Vert u_{\lambda}\\
 & \leq\max\left(\left\Vert \pi(h)\right\Vert ,(1-\delta)\left\Vert h\right\Vert \right).\end{align*}
\end{proof}
\begin{lem}
\label{lem:quasiCentralLemma} Suppose $I\vartriangleleft A.$ For
any approximate unit $u_{\lambda}$ of $I$ quasicentral for $A,$
any $a$ in $A,$ and any real $0\leq\delta\leq1,$\[
\limsup_{\lambda}\left\Vert a\left(1-\delta u_{\lambda}\right)^{\frac{1}{2}}\right\Vert \leq\max\left(\left\Vert \pi(a)\right\Vert ,\left(1-\delta\right)^{\frac{1}{2}}\left\Vert a\right\Vert \right).\]
\end{lem}
\begin{proof}
Using the quasicentral property and Lemma~\ref{lem:approxUnitsLemma},\begin{align*}
\limsup_{\lambda}\left\Vert a\left(1-\delta u_{\lambda}\right)^{\frac{1}{2}}\right\Vert ^{2} & =\limsup_{\lambda}\left\Vert a^{*}a\left(1-\delta u_{\lambda}\right)\right\Vert \\
 & =\limsup_{\lambda}\left\Vert a^{*}a\left(1-u_{\lambda}\right)+(1-\delta)a^{*}au_{\lambda}\right\Vert \\
 & =\limsup_{\lambda}\left\Vert \left(1-u_{\lambda}\right)^{\frac{1}{2}}a^{*}a\left(1-u_{\lambda}\right)^{\frac{1}{2}}+(1-\delta)u_{\lambda}^{\frac{1}{2}}a^{*}au_{\lambda}^{\frac{1}{2}}\right\Vert \\
 & \leq\max\left(\left\Vert \pi(a^{*}a)\right\Vert ,(1-\delta)\left\Vert a^{*}a\right\Vert \right)\\
 & =\left(\max\left(\left\Vert \pi(a)\right\Vert ,(1-\delta)\left\Vert a\right\Vert \right)\right)^{2}.\end{align*}
\end{proof}
\begin{thm}
\label{thm:normTechnicalLemma} Suppose $I\vartriangleleft A,$ that
$u_{\lambda}$ is a approximate unit $u_{\lambda}$ for $I$ quasicentral
for $A,$ and $a$ is in $A.$ If $f$ is a continuous function of
$[0,1]$ so that \[
1=f(0)\geq f(t)\geq f(1)\geq0\]
then\[
\limsup_{\lambda}\left\Vert af\left(u_{\lambda}\right)\right\Vert \leq\max\left(\left\Vert \pi(a)\right\Vert ,f\left(1\right)\left\Vert a\right\Vert \right).\]
\end{thm}
\begin{proof}
Let $\delta=1-f\left(1\right)^{2}$ and\[
g(t)=\delta^{-1}\left(1-f\left(t\right)^{2}\right).\]
This function is continuous and\[
0=g(0)\leq g(t)\leq g(1)=1\]
so $g\left(u_{\lambda}\right)$ is also a quasicentral approximate
unit. By the Lemma~\ref{lem:quasiCentralLemma},\begin{align*}
\limsup_{\lambda}\left\Vert af\left(u_{\lambda}\right)\right\Vert  & =\limsup_{\lambda}\left\Vert a\left(1-\delta g\left(u_{\lambda}\right)\right)^{\frac{1}{2}}\right\Vert \\
 & \leq\max\left(\left\Vert \pi(a)\right\Vert ,\left(1-\delta\right)^{\frac{1}{2}}\left\Vert a\right\Vert \right)\\
 & =\max\left(\left\Vert \pi(a)\right\Vert ,f\left(1\right)\left\Vert a\right\Vert \right).\end{align*}

\end{proof}

\section{Lifting Softened Homogeneous Relations\label{sec:Lifting-Softened-Homogeneous}}

We will consider $*$-polynomials in infinitely many variables that
are homogeneous in some finite subset of the variables. These we take
to be the first $r$-variables, which we label $x_{1},\ldots,x_{r},$
and the remaining variables we label $y_{1},y_{2},\ldots.$ We also
use the $n$-tuple notation $\mathbf{x}=(x_{1},\ldots,x_{r})$ and
$\mathbf{y}=(y_{1},y_{2},\ldots)$ and with a NC $*$-polynomial $p$
we use the notation \[
p(\mathbf{x},\mathbf{y})=p(x_{1},\ldots,x_{r},y_{1},y_{2},\ldots).\]
For scalar $t$ we use\[
t\mathbf{x}=(tx_{1},\ldots,tx_{r}).\]

We will say $p$ is \emph{$d$-homogeneous in the first $r$ variables
}if\[
p(t\mathbf{x},\mathbf{y})=t^{d}p(\mathbf{x},\mathbf{y})\]
for all real scalars $t.$ In other words, in each monomial the $x_{j}$
and $x_{j}^{*}$ appear collectively $d$ times. As $d$ is not necessarily
the degree of $p$ we call $d$ the \emph{degree of homogeneity of
$p.$}

We will rather quickly run out of letters if we insist on other symbols
when evaluating $p$ on specific elements of a $C^{*}$-algebra $A.$
Given $x_{1},\ldots,x_{r}$ in $A$ and $m$ in $A$ or $\widetilde{A}$
we define \[
m\mathbf{x}=(mx_{1},\ldots,mx_{r}).\]
(We use $\widetilde{A}$ to denote the unitization of $A.$) If $\varphi:A\rightarrow B$
then \[
\varphi(\mathbf{x})=(\varphi(x_{1}),\ldots,\varphi(x_{r}))\]
and\[
\varphi(\mathbf{y})=(\varphi(y_{1}),\varphi(y_{2}),\ldots).\]
If $z_{1},\ldots,z_{r}$ are in $A$ then $\mathbf{z}\leq\mathbf{y}$
shall mean $z_{j}\leq x_{j}$ for $j=1,\ldots,r.$ If $\mathbf{0}\leq\mathbf{x}$
then $\mathbf{x}^{\frac{1}{2}}$ shall denote $\left(x_{1}^{\frac{1}{2}},\ldots,x_{r}^{\frac{1}{2}}\right).$
For $a_{1},\ldots,a_{n}$ and $b_{1},\ldots,b_{n}$ in $A$ we use
the notation \[
a\bullet b=a_{1}b_{1}+\cdots+a_{n}b_{n}.\]

Recall $\pi$ is our generic notation for the quotient map $A\rightarrow A/I.$
\begin{thm}
\label{thm:MainTheorem}Suppose $p_{1},\ldots,p_{J}$ are NC $*$-polynomials
in infinitely many variables that are homogeneous in the first $r$
variables, with each degree of homogeneity $d_{j}$ at least one.
Suppose $C_{j}>0$ are constants. For every $C^{*}$-algebra $A$
and $I\vartriangleleft A$ an ideal, given $x_{1},\ldots,x_{r}$ and
$y_{1},y_{2},\ldots$ in $A$ with $\mathbf{0}\leq\mathbf{x}$ and
\[
\left\Vert p_{j}\left(\pi(\mathbf{x}),\pi(\mathbf{y})\right)\right\Vert \leq C_{j},\]
there are $z_{1},\ldots z_{r}$ in $A$ with $\mathbf{0}\leq\mathbf{z}\leq\mathbf{x}$
and $\pi(\mathbf{z})=\pi(\mathbf{x})$ and\begin{equation}
\left\Vert p_{j}\left(\mathbf{z},\mathbf{y}\right)\right\Vert \leq C_{j}.\label{eq:liftedSoftRelation}\end{equation}
\end{thm}
\begin{proof}
Our proof is modeled on that from \cite{AkePederIdealPerturb}. 

We start by performing the easier lifting where $C_{j}$ in (\ref{eq:liftedSoftRelation})
is replaced by $\left(1+\epsilon_{1}\right)C_{j}.$ We pick $\epsilon_{1}$
later, but it will be positive. Since $C_{j}$ is not allowed to be
zero, $\left(1+\epsilon_{1}\right)C_{j}$ will be strictly larger
than $C_{j}.$

Let $u_{\lambda}$ be any approximate unit $u_{\lambda}$ for $I$
that is quasicentral for $A.$ By quasicentrality and the homogeneity
in $\mathbf{x},$ we have \begin{align*}
\lim_{\lambda}\left\Vert p_{j}\left(\mathbf{x}^{\frac{1}{2}}\bullet\left(1-u_{\lambda}\right)\mathbf{x}^{\frac{1}{2}},\mathbf{y}\right)\right\Vert  & =\lim_{\lambda}\left\Vert p_{j}\left(\mathbf{x},\mathbf{y}\right)\left(1-u_{\lambda}\right)^{d_{j}}\right\Vert \\
 & \leq\lim_{\lambda}\left\Vert p_{j}\left(\mathbf{x},\mathbf{y}\right)\left(1-u_{\lambda}\right)\right\Vert \\
 & =\left\Vert \pi\left(p_{j}\left(\mathbf{x},\mathbf{y}\right)\right)\right\Vert \\
 & =\left\Vert p_{j}\left(\pi\left(\mathbf{x}\right),\pi\left(\mathbf{y}\right)\right)\right\Vert \\
 & \leq C_{j}.\end{align*}
We define \[
z_{k}^{(1)}=x_{k}^{\frac{1}{2}}\left(1-u_{\lambda_{1}}\right)x_{k}^{\frac{1}{2}}\]
where $\lambda_{1}$ is large enough to give us\[
\left\Vert p_{j}\left(\mathbf{z}^{(1)},\mathbf{y}\right)\right\Vert \leq\left(1+\epsilon_{1}\right)C_{j}\]
for $j=1,\ldots,J.$ Clearly $0\leq\mathbf{z}^{(1)}\leq\mathbf{x}$
and $\pi\left(\mathbf{z}^{(1)}\right)=\pi\left(\mathbf{x}\right).$

We will create ever better lifts by defining\[
z_{k}^{(2)}=\left(z_{k}^{(1)}\right)^{\frac{1}{2}}\left(1-\delta_{2}u_{\lambda_{2}}\right)\left(z_{k}^{(1)}\right)^{\frac{1}{2}},\]
and so forth. For consistency, we let $\delta_{1}=1$ and $\mathbf{z}^{(0)}=\mathbf{x}.$
We choose $\delta_{1}>\delta_{2}>\dots$ all positive with \begin{equation}
\sum_{c=1}^{\infty}\delta_{c}<\infty.\label{eq:deltasSummable}\end{equation}
We set $\epsilon_{1}$ and the rest of a sequence $\epsilon_{c}$
by requiring\[
\left(1-\delta_{c+1}\right)\left(1+\epsilon_{c}\right)=1.\]
Notice the $\epsilon_{c}$ are positive and decreasing to zero.

Assume we have found $\mathbf{z}^{(1)}$ through $\mathbf{z}^{(c-1)}$
with \begin{equation}
\mathbf{0}\leq\mathbf{z}^{(c-1)}\leq\mathbf{z}^{(c-2)}\cdots\leq\mathbf{z}^{(1)}\leq\mathbf{x}\label{eq:orderOfTheZ}\end{equation}
\begin{equation}
\pi\left(\mathbf{z}^{(c-1)}\right)=\pi\left(\mathbf{z}^{(c-2)}\right)=\cdots=\pi\left(\mathbf{z}^{(1)}\right)=\pi\left(\mathbf{x}\right)\label{eq:LIftCorrectThings}\end{equation}
and\begin{equation}
\left\Vert p_{j}\left(\mathbf{z}^{(w)},\mathbf{y}\right)\right\Vert \leq\left(1+\epsilon_{w}\right)C_{j},\quad\left(w=1,\ldots,c-1\right).\label{eq:HowGoodTheNorms}\end{equation}
Moreover, assume the $\mathbf{z}^{(w)}$ have been constructed via
the formula \begin{equation}
\mathbf{z}^{(w)}=\left(\mathbf{z}^{(w-1)}\right)^{\frac{1}{2}}\left(1-\delta_{w}u_{\lambda_{w}}\right)\left(\mathbf{z}^{(w-1)}\right)^{\frac{1}{2}},\quad\left(w=1,\ldots,c-1\right).\label{eq:SandwichConstruction}\end{equation}
 Theorem~\ref{thm:normTechnicalLemma} tells us\begin{align*}
 & \limsup_{\lambda}\left\Vert p_{j}\left(\left(\mathbf{z}^{(c-1)}\right)^{\frac{1}{2}}\bullet\left(1-\delta_{c}u_{\lambda}\right)\left(\mathbf{z}^{(c-1)}\right)^{\frac{1}{2}},\mathbf{y}\right)\right\Vert \\
 & \hspace{1in}=\limsup_{\lambda}\left\Vert p_{j}\left(\mathbf{z}^{(c-1)},\mathbf{y}\right)\left(1-\delta_{c}u_{\lambda}\right)^{d_{j}}\right\Vert \\
 & \hspace{1in}\leq\limsup_{\lambda}\left\Vert p_{j}\left(\mathbf{z}^{(c-1)},\mathbf{y}\right)\left(1-\delta_{c}u_{\lambda}\right)\right\Vert \\
 & \hspace{1in}\leq\max\left(\left\Vert \pi\left(p_{j}\left(\mathbf{z}^{(c-1)},\mathbf{y}\right)\right)\right\Vert ,\left(1-\delta_{c}\right)\left\Vert p_{j}\left(\mathbf{z}^{(c-1)},\mathbf{y}\right)\right\Vert \right)\\
 & \hspace{1in}=\max\left(\left\Vert p_{j}\left(\pi\left(\mathbf{x}\right),\pi\left(\mathbf{y}\right)\right)\right\Vert ,\left(1-\delta_{c}\right)\left\Vert p_{j}\left(\mathbf{z}^{(c-1)},\mathbf{y}\right)\right\Vert \right)\\
 & \hspace{1in}\leq\max\left(C_{j},\left(1-\delta_{c}\right)\left(\left(1+\epsilon_{c-1}\right)C_{j}\right)\right)\\
 & \hspace{1in}=C_{j}\end{align*}
so we may choose $\lambda_{c}$ with\[
\left\Vert p_{j}\left(\left(\mathbf{z}^{(c-1)}\right)^{\frac{1}{2}}\bullet\left(1-\delta_{c}u_{\lambda_{c}}\right)\left(\mathbf{z}^{(c-1)}\right)^{\frac{1}{2}},\mathbf{y}\right)\right\Vert \leq\left(1+\epsilon_{c}\right)C_{j}.\]
We set\[
\mathbf{z}^{(c)}=\left(\mathbf{z}^{(c-1)}\right)^{\frac{1}{2}}\left(1-\delta_{c}u_{\lambda_{c}}\right)\left(\mathbf{z}^{(c-1)}\right)^{\frac{1}{2}}\]
 and the construction continues.

We wish to set ${\displaystyle z_{k}=\lim_{c}z_{j}^{(c)},}$ and we
may because\begin{align*}
\left\Vert z_{k}^{(c)}-z_{k}^{(c-1)}\right\Vert  & =\left\Vert \left(z_{k}^{(c-1)}\right)^{\frac{1}{2}}\left(\delta_{c}u_{\lambda_{c}}\right)\left(z_{k}^{(c-1)}\right)^{\frac{1}{2}}\right\Vert \\
 & \leq\left\Vert z_{k}^{(c-1)}\right\Vert \left\Vert \delta_{c}u_{\lambda_{c}}\right\Vert \\
 & \leq2\left\Vert x_{k}\right\Vert \delta_{c}.\end{align*}
Equations (\ref{eq:orderOfTheZ}), (\ref{eq:LIftCorrectThings}) and
(\ref{eq:HowGoodTheNorms}) give us $\mathbf{0}\leq\mathbf{z}\leq\mathbf{x}$
and $\pi(\mathbf{z})=\mathbf{x}$ and finally the norm conditions\[
\left\Vert p_{j}\left(\mathbf{z},\mathbf{y}\right)\right\Vert =\lim_{c}\left\Vert p_{j}\left(\mathbf{z}^{(c)},\mathbf{y}\right)\right\Vert \leq C_{j}.\]

\end{proof}
If we have soft $*$-polynomial relations involving self-adjoint variables
we can replace each by two positive variables. A variable that is
a contraction can be replaced by four positive variables. These replacements
will preserve any homogeneity in a subset of the variables. Thus we
can have a more flexible version of Theorem~\ref{thm:MainTheorem}.
As stated, Theorem~\ref{thm:MixedVariables} it is not a corollary
as we are very specific in how the lifts are adjusted for the different
types of variables.
\begin{thm}
\label{thm:MixedVariables}Suppose $p_{1},\ldots,p_{J}$ are NC $*$-polynomials
in infinitely many variables that are homogeneous in the first $r$
variables, with each degree of homogeneity $d_{j}$ at least one.
Suppose $C_{j}>0$ are constants. Suppose $S=\{1,\ldots,r\}$ is partitioned
as\[
S=S_{+}\cup S_{\mathrm{h}}\cup S_{\mathrm{g}},\]
we have positive constants $C_{j},$ nonnegative constants $D_{k},$
$E_{k},$ $F_{k}$ and $G_{k},$ and consider the relations\begin{equation}
0\leq x_{k}\leq D_{k}\quad(k\in S_{+})\label{eq:PositiveBound}\end{equation}
\begin{equation}
E_{k}\leq x_{k}\leq F_{k}\quad(k\in S_{\mathrm{h}})\label{eq:HermitianBound}\end{equation}
\begin{equation}
\left\Vert x_{k}\right\Vert \leq G_{k}\quad(k\in S_{\mathrm{g}})\label{eq:NormBound}\end{equation}
\begin{equation}
\left\Vert p_{j}\left(\mathbf{x},\mathbf{y}\right)\right\Vert \leq C_{j}.\label{eq:PolyNormBound}\end{equation}
For every $C^{*}$-algebra $A$ and $I\vartriangleleft A$ an ideal,
given $x_{1},\ldots,x_{r}$ and $y_{1},y_{2},\ldots$ in $A$ so that
$\left(\mathbf{x},\mathbf{y}\right)$ satisfies (\ref{eq:PositiveBound}-\ref{eq:NormBound})
and $\left(\pi\left(\mathbf{x}\right),\pi\left(\mathbf{y}\right)\right)$
satisfies (\ref{eq:PolyNormBound}), there are elements $z_{1},\ldots,z_{r}$
in $A$ so that $\left(\mathbf{z},\mathbf{y}\right)$ satisfy (\ref{eq:PositiveBound}-\ref{eq:PolyNormBound})
and $\pi\left(\mathbf{z}\right)=\pi\left(\mathbf{x}\right).$ Moreover,
it is possible to do so with \begin{align*}
z_{k} & =\left(x_{k}\right)^{\frac{1}{2}}m^{2}\left(x_{k}\right)^{\frac{1}{2}},\quad(k\in S_{+})\\
z_{k} & =mx_{k}m,\quad(k\in S_{h})\\
z_{k} & =x_{k}m^{2},\quad(k\in S_{\mathrm{g}})\end{align*}
 for some $m$ in $1+I$ with $0\leq m\leq1.$\end{thm}
\begin{proof}
Let $\epsilon_{c},$ $\delta_{c}$ and $u_{\lambda}$ be as before.
We modify the construction used for Theorem~\ref{thm:MainTheorem}
by requiring $m_{0}=1$ and \[
m_{c}=\left(1-\delta_{c}u_{\lambda_{c}}\right)m_{c-1}\left(1-\delta_{c}u_{\lambda_{c}}\right)\]
and \begin{align*}
z_{k}^{(c)} & =x_{k}^{\frac{1}{2}}m_{c}^{2}x_{k}^{\frac{1}{2}}\quad(k\in S_{+})\\
z_{k}^{(c)} & =m_{c}x_{k}m_{c}\quad(k\in S_{h})\\
z_{k}^{(c)} & =x_{k}m_{c}^{2}\quad(k\in S_{\mathrm{g}}).\end{align*}
We want\[
\left\Vert p_{j}\left(\mathbf{z}^{(c)},\mathbf{y}\right)\right\Vert \leq\left(1+\epsilon_{c}\right)C_{j}\]
given that we already have defined $m_{c-1}$ in $1+I$ to established\begin{equation}
\left\Vert p_{j}\left(\mathbf{z}^{(c-1)},\mathbf{y}\right)\right\Vert \leq\left(1+\epsilon_{c-1}\right)C_{j}.\label{eq:baseBound}\end{equation}
To unify the initial step and subsequent steps, we take $\epsilon_{0}$
large enough to force (\ref{eq:baseBound}) when $c=1.$ We need to
find the right $\lambda_{c}$ to define $z_{j}^{(c)}=w_{j}^{(\lambda_{c})}$
where\begin{align*}
w_{k}^{(\lambda)} & =x_{k}^{\frac{1}{2}}\left(\left(1-\delta_{c}u_{\lambda}\right)m_{c-1}\left(1-\delta_{c}u_{\lambda}\right)\right)^{2}x_{k}^{\frac{1}{2}}\quad(k\in S_{+})\\
w_{k}^{(\lambda)} & =\left(\left(1-\delta_{c}u_{\lambda}\right)m_{c-1}\left(1-\delta_{c}u_{\lambda}\right)\right)x_{k}\left(\left(1-\delta_{c}u_{\lambda}\right)m_{c-1}\left(1-\delta_{c}u_{\lambda}\right)\right)\quad(k\in S_{h})\\
w_{k}^{(\lambda)} & =x_{k}\left(\left(1-\delta_{c}u_{\lambda}\right)m_{c-1}\left(1-\delta_{c}u_{\lambda}\right)\right)^{2}\quad(k\in S_{\mathrm{g}})\end{align*}
and see\[
\lim_{\lambda}\left\Vert w_{j}^{(\lambda)}-\left(1-\delta_{c}u_{\lambda}\right)^{4}z_{j}^{(c-1)}\right\Vert =0.\]
Therefore \begin{align*}
\limsup_{\lambda}\left\Vert p_{j}\left(\mathbf{w}^{(\lambda)},\mathbf{y}\right)\right\Vert  & =\limsup_{\lambda}\left\Vert p_{j}\left(\left(1-\delta_{c}u_{\lambda}\right)^{4}\mathbf{z}^{(c-1)},\mathbf{y}\right)\right\Vert \\
 & =\limsup_{\lambda}\left\Vert p_{j}\left(\mathbf{z}^{(c-1)},\mathbf{y}\right)\left(1-\delta_{c}u_{\lambda}\right)^{4d_{j}}\right\Vert \\
 & \leq\limsup_{\lambda}\left\Vert p_{j}\left(\mathbf{z}^{(c-1)},\mathbf{y}\right)\left(1-\delta_{c}u_{\lambda}\right)\right\Vert \\
 & \leq\max\left(\left\Vert p_{j}\left(\pi\left(\mathbf{x}\right),\pi\left(\mathbf{y}\right)\right)\right\Vert ,\left(1-\delta_{c}\right)\left\Vert p_{j}\left(\mathbf{z}^{(c-1)},\mathbf{y}\right)\right\Vert \right)\\
 & \leq C_{j}\end{align*}
 and it is possible to chose the needed $\lambda_{c}.$ 

It is clear that $m_{c}$ stays in $i+I,$ so these are all lifts
of the original $\mathbf{x}.$ What is left to check is that $m={\displaystyle \lim_{c}m_{c}}$
exists. Indeed it does, as\begin{align*}
\left\Vert m_{c}-m_{c-1}\right\Vert  & =\left\Vert \left(1-\delta_{c}u_{\lambda_{c}}\right)m_{c-1}\left(1-\delta_{c}u_{\lambda_{c}}\right)-m_{c-1}\right\Vert \\
 & \leq\left\Vert \left(\delta_{c}u_{\lambda_{c}}\right)m_{c-1}\left(1-\delta_{c}u_{\lambda_{c}}\right)\right\Vert +\left\Vert m_{c-1}\left(\delta_{c}u_{\lambda_{c}}\right)\right\Vert \\
 & \leq2\left\Vert \delta_{c}u_{\lambda_{c}}m_{c-1}\right\Vert \\
 & \leq2\delta_{c}.\end{align*}

\end{proof}
We get from Theorem~\ref{thm:MixedVariables} a myriad of projective
$C^{*}$-algebras, simply by adding relations such as $-1\leq y_{j}\leq1$
that are liftable and that impose a norm restriction forcing the universal
$C^{*}$-algebra to exist. We generally add the relation $\left\Vert y_{j}\right\Vert =0$
to most of the $y_{j}$ so as to be working with a finitely generated
projective $C^{*}$-algebra.

\section{Soft Versions of Known Projectives\label{sec:Soft-Versions-of}}

Consider $C_{0}(X)$ where $X^{+}$ is a (finite) tree. The presentation
in \cite{LoringProjectiveCstar} for the projective $C^{*}$-algebra
$C_{0}(X)$ was based on a partial order $\preceq$ on $\{1,\ldots,s\}.$
This was not a general partial order, it had to be the partial order
on the non-root vertices determined by paths away from the root. Let
us call such a relation a tree order.

The presentation associated to the tree order $\preceq$ had generators
$h_{1},\ldots,h_{s}$ and relations\begin{align*}
0\leq h_{j} & \leq1,\quad(j=1,\ldots,s)\\
h_{i}h_{j} & =h_{j},\quad(\mbox{if }i\prec j)\\
h_{i}h_{j} & =0,\quad(\mbox{if }i\not\preceq j\mbox{ and }j\not\preceq i).\end{align*}
The last two lines of relations are not generally homogeneous in any
subset of the variables. Taking advantage of the precise way the liftings
are modified in our main theorem, we can still lift the soft version
of these.

The concrete function in $C_{0}(X)$ that corresonds to the abstract
$h_{j}$ only varies along one egdge, the edge incident to the $j$th
vertex on a path from that vertex to the root. From this point of
view, it makes more sense to index the generators by the edges (as
in \cite{LoringProjectiveCstar}) but in graph theory, tree orders
are on the vertices.

We do need to make two changes to the relations. Some are redundant,
since for positive elements \[
h_{1}h_{2}=h_{2}\ \&\ h_{3}h_{1}=0\implies h_{3}h_{2}=0\]
 and\[
h_{1}h_{2}=h_{2}\ \&\ h_{3}h_{1}=h_{1}\implies h_{3}h_{2}=h_{2}.\]
Let us swich to the indexing being over the non-root vertices. We
then can speak of $i$ being a child of $j,$ meaning $j\preceq i$
and \[
j\preceq k\preceq i\implies k=j\mbox{ or }k=i.\]
The only relations we need are those that ask that the parent act
as a unit on the child and that two children of the same parent must
be orthogonal. We call two children of the same vertex siblings, of
course. Children of children, and so forth, we call descendents.

A second change is we replace $h_{i}h_{j}=h_{j}$ by\[
\left(h_{i}-1\right)h_{j}\left(h_{i}-1\right)=0.\]
 
\begin{thm}
Suppose $\preceq$ is a tree order on $\{1,\ldots,s\}.$ The relations
\begin{align*}
0\leq h_{j} & \leq1,\quad(j=1,\ldots,s)\\
\left\Vert \left(h_{i}-1\right)h_{j}\left(h_{i}-1\right)\right\Vert  & \leq\epsilon,\quad(\mbox{if }j\mbox{ is a child of }i)\\
\left\Vert h_{i}h_{j}\right\Vert  & \leq\epsilon,\quad(\mbox{if }i\mbox{ and }j\mbox{ are siblings})\end{align*}
 are liftable.
\end{thm}
The following, slightly stronger result is more easily proven.
\begin{thm}
\label{thm:SoftTrees}Suppose $\preceq$ is a tree order on $\{1,\ldots,s\}$
and $\epsilon>0.$ For every $C^{*}$-algebra $A$ and $I\vartriangleleft A$
an ideal, given $h_{1},\ldots,h_{s}$ in $A$ with $\mathbf{0}\leq\mathbf{h}\leq\mathbf{1}$
and \begin{align*}
\left\Vert \pi\left(\left(h_{i}-1\right)h_{j}\left(h_{i}-1\right)\right)\right\Vert  & \leq\epsilon,\quad(\mbox{if }j\mbox{ is a child of }i)\\
\left\Vert \pi\left(h_{i}\right)\pi\left(h_{j}\right)\right\Vert  & \leq\epsilon,\quad(\mbox{if }i\mbox{ and }j\mbox{ are siblings})\end{align*}
there are $k_{1},\ldots k_{s}$ in $A$ with $\mathbf{0}\leq\mathbf{k}\leq\mathbf{h}$
and $\pi(\mathbf{k})=\pi(\mathbf{h})$ and\begin{align*}
\left\Vert \left(k_{i}-1\right)k_{j}\left(k_{i}-1\right)\right\Vert  & \leq\epsilon,\quad(\mbox{if }j\mbox{ is a child of }i)\\
\left\Vert k_{i}k_{j}\right\Vert  & \leq\epsilon,\quad(\mbox{if }i\mbox{ and }j\mbox{ are siblings}).\end{align*}
\end{thm}
\begin{proof}
The very trivial base case for our proof by induction is the case
of zero generators.

Re-indexing, we may assume the minimal elements (i.e.~elders or vertices
closest to the root) in this partial order are $\{1,\ldots,r\}.$
If $m$ and $n$ are descendents of different minimal elements $i$
and $j$ then there are no relations involving both $h_{m}$ and $h_{n}.$
The relations not involving the minimal elements are a disjoint union
of relations of the type in the statement of the theorem.

The relations we need that involve the minimal elements are $0\leq h_{j}\leq1$
for $j=1,\ldots,r$ and\begin{align*}
\left\Vert \left(h_{i}-1\right)h_{m}\left(h_{i}-1\right)\right\Vert  & \leq\epsilon,\quad(m\mbox{ is the child of }i)\\
\left\Vert k_{i}k_{j}\right\Vert  & \leq\epsilon,\quad(1\leq i<j\leq r).\end{align*}
These are homogeneous in $\{h_{1},\ldots,h_{r}\}.$ By Theorem~\ref{thm:MainTheorem}
there are $k_{1},\ldots,k_{r}$ in $A$ with $0\leq k_{j}\leq h_{j}$
and $\pi(h_{j})=\pi(h_{j})$ for $j\leq r$ and \begin{align*}
\left\Vert \left(k_{i}-1\right)h_{m}\left(k_{i}-1\right)\right\Vert  & \leq\epsilon,\quad(m\mbox{ is the child of }i)\\
\left\Vert k_{i}k_{j}\right\Vert  & \leq\epsilon,\quad(1\leq i<j\leq r).\end{align*}
The induction hypothesis tells us there are $k_{r+1},\ldots,k_{s}$
with $0\leq k_{m}\leq h_{m}$ and $\pi(k_{m})=\pi(h_{m})$ for $m>r$
with all the relations not involving indices $\{1,\ldots,r\}.$ We
might have lost the relations between some $h_{i}$ and $h_{m}$ with
$m$ a child of $i,$ but we have not, since\[
\left(h_{i}-1\right)k_{m}\left(h_{i}-1\right)\leq\left(h_{i}-1\right)h_{m}\left(h_{i}-1\right).\]

\end{proof}
A rather different example, with a similar proof, is a soft version
of the projective $C^{*}$-algebra\[
C^{*}\left\langle h,k,x\left|\begin{array}{c}
hk=0\\
0\leq\left[\begin{array}{cc}
\mathbb{1}-h & x^{*}\\
x & k\end{array}\right]\leq1\end{array}\right.\right\rangle \]
considered in \cite{LoringProjectiveKtheory}. (For a detailed explanation
of how the second relation is valid, see \cite{LoringCstarRelations}.)
\begin{thm}
\label{thm:projectiveLike_qC}For any positive $\epsilon,$ the $C^{*}$-algebra\[
C^{*}\left\langle h,k,x\left|\begin{array}{c}
\left\Vert hk\right\Vert \leq\epsilon\\
0\leq\left[\begin{array}{cc}
\mathbb{1}-h & x^{*}\\
x & k\end{array}\right]\leq1\end{array}\right.\right\rangle \]
is projective.\end{thm}
\begin{proof}
Suppose $h,$ $k$ and $x$ are in $A,$ which we may assume is unital,
are such that\[
\left\Vert \pi\left(h\right)\pi\left(k\right)\right\Vert \leq\epsilon\]
and\[
0\leq\left[\begin{array}{cc}
1-\pi\left(h\right) & \pi\left(x\right)^{*}\\
\pi\left(x\right) & \pi\left(k\right)\end{array}\right]\leq1.\]
 We know positive contractions lift to positive contractions from
$\mathbf{M}_{2}(A/I)$ to $\mathbf{M}_{2}(A)$ and so we can find
$\hat{h},$ $\hat{k}$ and $\hat{x}$ in $A$ so that $\pi\left(\hat{h}\right)=\pi\left(h\right),$
$\pi\left(\hat{k}\right)=\pi\left(k\right),$ $\pi\left(\hat{x}\right)=\pi\left(x\right)$
and \[
0\leq\left[\begin{array}{cc}
1-\hat{h} & \hat{x}^{*}\\
\hat{x} & \hat{k}\end{array}\right]\leq1.\]
 The polynomial $hk$ is homogeneous in $k$ so Theorem~\ref{thm:MixedVariables}
tells us there is $0\leq m\leq1$ in $1+I$ so that $\left\Vert \hat{h}m\hat{k}m\right\Vert \leq\epsilon.$
Let $\bar{h}=\hat{h},$ $\bar{x}=m\hat{x}$ and $\bar{k}=m\hat{k}m.$
These are still lifts of $h,$ $x$ and $k,$ and now $\left\Vert \bar{h}\bar{k}\right\Vert \leq\epsilon$
and\[
\left[\begin{array}{cc}
1-\bar{h} & \bar{x}^{*}\\
\bar{h} & \bar{h}\end{array}\right]=\left[\begin{array}{cc}
1 & 0\\
0 & m\end{array}\right]\left[\begin{array}{cc}
1-\hat{h} & \hat{x}^{*}\\
\hat{x} & \hat{k}\end{array}\right]\left[\begin{array}{cc}
1 & 0\\
0 & m\end{array}\right]\]
implies\[
0\leq\left[\begin{array}{cc}
1-\bar{h} & \bar{x}^{*}\\
\bar{h} & \bar{h}\end{array}\right]\leq1.\]

\end{proof}

\section{Fattened Curves in Various NC unit balls \label{sec:Fattened-Curves-in}}
\begin{thm}
\label{thm:setsInSquare} Suppose $p_{1},\ldots,p_{J}$ are NC $*$-polynomials
in $x_{1},\ldots,x_{s}.$ Suppose $1\leq r\leq s$ and each $p_{j}$
is homogeneous in $\{x_{1},\ldots,x_{r}\}$ with degree of homogeneity
$d_{j}\geq1.$ For $\epsilon>0,$ the $C^{*}$-algebra\emph{\[
A_{\epsilon}=C^{*}\left\langle x_{1},\ldots,x_{s}\left|\,\begin{array}{c}
\left\Vert x_{k}\right\Vert \leq1,\, k=1,\ldots,s\\
\left\Vert p_{j}\left(x_{1},\ldots,x_{s}\right)\right\Vert \leq\epsilon,\, j=1,\ldots,J\end{array}\right.\right\rangle \]
 is projective.}\end{thm}
\begin{proof}
This is immediate consequence of Theorem~\ref{thm:MixedVariables}.
\end{proof}
For a single NC $*$-polynomial $p,$ we can think of \emph{\begin{equation}
C^{*}\left\langle x_{1},\ldots,x_{s}\left|\,\begin{array}{c}
\left\Vert x_{k}\right\Vert \leq1\\
\left\Vert p\left(x_{1},\ldots,x_{s}\right)\right\Vert \leq\epsilon\end{array}\right.\right\rangle \label{eq:curveInSquare}\end{equation}
}as a approximate zero locus of a NC curve intersected with the NC
unit square. Likewise\emph{ }we can think of\begin{equation}
C^{*}\left\langle x_{1},\ldots,x_{s}\left|\,\begin{array}{c}
\left\Vert {\displaystyle \sum_{k=1}^{s}x_{k}^{*}x_{k}}\right\Vert \leq1\\
\left\Vert p\left(x_{1},\ldots,x_{s}\right)\right\Vert \leq\epsilon\end{array}\right.\right\rangle \label{eq:curveInBall}\end{equation}
as a approximate zero locus of a NC curve with the NC unit ball. Notice
that the ``row contraction'' condition \emph{$\left\Vert {\displaystyle \sum x_{k}^{*}x_{k}}\right\Vert \leq1$
}implies $\left\Vert x_{k}\right\Vert \leq1$ so we can still apply
Theorem~\ref{thm:setsInSquare}.\emph{ }For $\epsilon>0,$ and with
$p$ homogeneous in $x_{1},\ldots,x_{s}$ we find (\ref{eq:curveInSquare})
and (\ref{eq:curveInBall}) define projective $C^{*}$-agebras.

We will see that it is possible to work with other unit balls, not
just the ones corresponding to the $\ell^{2}$ and $\ell^{\infty}$
norms.
\begin{lem}
\label{lem:approxRoots} Suppose $0<\alpha<\infty$ is a scalar. For
every $\epsilon>0$ there is a $\delta>0$ so that for any two positive
contractions in any $C^{*}$-algebra, \[
\left\Vert hk-kh\right\Vert \leq\delta\implies\left\Vert \left(hkh\right)^{\alpha}-k^{\alpha}h^{2\alpha}\right\Vert \leq\epsilon.\]
\end{lem}
\begin{proof}
This can be rephrased so it becomes a special case of Lemma~10 of
\cite{LoringProjectiveKtheory}, but it is easier to just revise the
proof. We know for nonnegative scalars $\left(xyx\right)^{\alpha}=y^{\alpha}x^{2\alpha}$
so by spectral theory, \[
hk=kh\implies\left(hkh\right)^{\alpha}=k^{\alpha}h^{2\alpha}.\]
If the lemma is false, there must be some $\epsilon_{0}$ and $h_{n}$
and $k_{n}$ in $A_{n}$ with $0\leq h_{n}\leq1$ and $0\leq k_{n}\leq1$
and \[
\left\Vert h_{n}k_{n}-k_{n}h_{n}\right\Vert \leq\frac{1}{n}\]
 and \[
\left\Vert \left(h_{n}k_{n}h_{n}\right)^{\alpha}-k_{n}^{\alpha}h_{n}^{2\alpha}\right\Vert \geq\epsilon_{0}.\]
This creates an element in\[
\left.\prod A_{n}\right/\bigoplus A_{n}\]
with $hk=kh$ and $\left(hkh\right)^{\alpha}-k^{\alpha}h^{2\alpha}\neq0,$
a contradiction.\end{proof}
\begin{thm}
\label{thm:AltNCballs} Suppose $r$ is a natural number. For $0<p<\infty$
define \[
B_{p}=C^{*}\left\langle x_{1},\dots,x_{r}\left|\strut\,\left\Vert \sum_{k=1}^{r}\left(x_{k}x_{k}^{*}\right)^{\frac{p}{2}}\right\Vert \leq1\right.\right\rangle \]
If $0<p\leq\infty$ then $B_{p}$ is projective. \end{thm}
\begin{proof}
Suppose $I\vartriangleleft A$ with quasicentral approximate unit
$u_{\lambda}.$ Suppose $x_{k}$ are in $A$ with $\pi\left(x_{k}\right)$
in $A/I$ satisfying\[
\left\Vert \sum_{k=1}^{r}\left(\pi\left(x_{k}\right)\left(\pi\left(x_{k}\right)\right)^{*}\right)^{\frac{p}{2}}\right\Vert \leq1.\]
Let $z_{k}^{(0)}=x_{k}$ and let $\epsilon_{0}$ be sufficiently large
so as to have\[
\left\Vert \sum_{k=1}^{r}\left(z_{k}^{(1)}\left(z_{k}^{(1)}\right)^{*}\right)^{\frac{p}{2}}\right\Vert \leq1+\epsilon_{0}.\]
Choose $\delta_{c}$ a positive sequence decreasing to zero with $\delta_{1}=1$
and\[
\sum_{c=1}^{\infty}\left(1-\left(1-\delta_{c}\right)^{\frac{1}{p}}\right)<\infty.\]
Define $\epsilon_{c}\searrow0$ for $c\geq1$ by the formula \[
\left(1-\delta_{c+1}\right)\left(1+\epsilon_{c}\right)=1.\]

Assume we have found $\mathbf{z}^{(c-1)}$ with\[
\pi\left(\mathbf{z}^{(c-1)}\right)=\pi\left(\mathbf{x}\right)\]
and\[
\left\Vert \sum_{k=1}^{r}\left(z_{k}^{(c-1)}\left(z_{k}^{(c-1)}\right)^{*}\right)^{\frac{p}{2}}\right\Vert \leq1+\epsilon_{c}.\]
Using Lemma~\ref{lem:approxRoots} and Theorem~\ref{thm:normTechnicalLemma}
we find \begin{align*}
 & \limsup_{\lambda}\left\Vert \sum_{k=1}^{r}\left(\left(1-\delta_{c}u_{\lambda}\right)^{\frac{1}{p}}z_{k}^{(c-1)}\left(z_{k}^{(c-1)}\right)^{*}\left(1-\delta_{c}u_{\lambda}\right)^{\frac{1}{p}}\right)^{\frac{p}{2}}\right\Vert \\
 & \hspace{1in}=\limsup_{\lambda}\left\Vert \sum_{k=1}^{r}\left(z_{k}^{(c-1)}\left(z_{k}^{(c-1)}\right)^{*}\right)^{\frac{p}{2}}\left(1-\delta_{c}u_{\lambda}\right)\right\Vert \\
 & \hspace{1in}\leq1.\end{align*}
We can chose $\lambda_{c}$ and set\[
z_{k}^{(c)}=\left(1-\delta_{c}u_{\lambda_{c}}\right)^{\frac{1}{p}}z_{k}^{(c-1)}\]
where $\lambda_{c}$ is large enough to ensure\[
\left\Vert \sum_{k=1}^{r}\left(z_{k}^{(c)}\left(z_{k}^{(c)}\right)^{*}\right)^{\frac{p}{2}}\right\Vert \leq1+\epsilon_{c}.\]
 The $z_{j}^{(c)}$ converge because \begin{align*}
\left\Vert z_{k}^{(c)}-z_{k}^{(c-1)}\right\Vert  & =\left\Vert \left(1-\delta_{c}u_{\lambda_{c}}\right)^{\frac{1}{p}}z_{k}^{(c-1)}-z_{k}^{(c-1)}\right\Vert \\
 & \leq\left\Vert \left(1-\delta_{c}u_{\lambda_{c}}\right)^{\frac{1}{p}}-1\right\Vert \left\Vert z_{k}^{(c-1)}\right\Vert \\
 & \leq\left(1-\left(1-\delta_{c}\right)^{\frac{1}{p}}\right)\left\Vert x_{k}\right\Vert .\end{align*}

\end{proof}
We still cannot tell if all the $B_{p}$ are isomorphic. They do interact
with homogeneous $*$-polynomials in about the same fashion as the
usual unit ball. The set of NC $*$-polynomials that we know we can
mix with the nonstandard unit ball condition depends on $p.$ We have
no idea if this is a limitation of our methods, or a real limitation.
\begin{thm}
\label{thm:setsInNCball-lower} Suppose $r\leq s$ and $p_{1},\ldots,p_{J}$
are NC $*$-polynomials in $x_{1},\ldots,x_{s},$ each homogeneous
in $\{x_{1},\ldots,x_{r}\}$ with degree of homogeneity $d_{j}$ at
least one. For $C_{j}>0$ and $0<q\leq2$ the $C^{*}$-algebra\emph{\[
A_{\epsilon}=C^{*}\left\langle x_{1},\ldots,x_{s}\left|\,\begin{array}{c}
\left\Vert {\displaystyle \sum_{k=1}^{s}\left(x_{k}x_{k}^{*}\right)^{\frac{q}{2}}}\right\Vert \leq1\\
\left\Vert p_{j}\left(x_{1},\ldots,x_{s}\right)\right\Vert \leq C_{j},\, j=1,\ldots,J\end{array}\right.\right\rangle \]
 is projective.}\end{thm}
\begin{proof}
Suppose we are given $\pi:A\rightarrow A/I$ with $x_{1},\ldots,x_{s}$
in $A$ with\[
\left\Vert \sum_{k=1}^{s}\left(\pi\left(x_{k}x_{k}^{*}\right)\right)^{\frac{q}{2}}\right\Vert \leq1\]
and\[
\left\Vert p_{j}\left(\pi\left(x_{1}\right),\ldots,\pi\left(x_{s}\right)\right)\right\Vert \leq C_{j}.\]
We first apply Theorem~\ref{thm:AltNCballs} to find $y_{1},\ldots,y_{s}$
in $A$ with $\pi(y_{k})=\pi(x_{k})$ and\[
\left\Vert \sum_{k=1}^{s}\left(y_{k}y_{k}^{*}\right)^{\frac{q}{2}}\right\Vert \leq1.\]
Theorem~\ref{thm:MixedVariables} gives us $z_{1},\ldots,z_{s}$
in $A$ with $\pi(z_{j})=\pi(x_{j})$ and\[
\left\Vert p_{j}\left(z_{1},\ldots,z_{s}\right)\right\Vert \leq C_{j},\]
but also with $z_{k}=y_{k}m$ for $k\leq r$ and $z_{k}=y_{k}$ for
$k>r,$ where $0\leq m\leq1.$ Therefore\[
z_{k}z_{k}^{*}=y_{k}m^{2}y_{k}^{*}\leq y_{k}y_{k}^{*}\]
for $1\leq k\leq r.$ Since for $q\le2$ the function $t^{q/2}$ is
operator-monotone we get \[
\sum_{k=1}^{s}\left(z_{k}z_{k}^{*}\right)^{\frac{q}{2}}\leq\sum_{k=1}^{s}\left(y_{k}y_{k}^{*}\right)^{\frac{q}{2}}\leq1.\]
\end{proof}
\begin{thm}
\label{thm:setsInNCball-higher} Suppose $p_{1},\ldots,p_{J}$ are
homogeneous, degree-$d_{j}$ NC $*$-polynomials in $x_{1},\ldots,x_{r}$
with $d_{j}\geq1.$ For $C_{j}>0$ and $2<q<\infty$ the $C^{*}$-algebra\emph{\[
A_{\epsilon}=C^{*}\left\langle x_{1},\ldots,x_{r}\left|\,\begin{array}{c}
\left\Vert {\displaystyle \sum_{k=1}^{r}\left(x_{k}x_{k}^{*}\right)^{\frac{q}{2}}}\right\Vert \leq1\\
\left\Vert p_{j}\left(x_{1},\ldots,x_{r}\right)\right\Vert \leq C_{j},\, j=1,\ldots,J\end{array}\right.\right\rangle \]
 is projective.}\end{thm}
\begin{proof}
Suppose we are given $\pi:A\rightarrow A/I$ with $x_{1},\ldots,x_{r}$
in $A$ with\[
\left\Vert \sum_{k=1}^{r}\left(\pi\left(x_{k}x_{k}^{*}\right)\right)^{\frac{q}{2}}\right\Vert \leq1\]
and\[
\left\Vert p_{j}\left(\pi\left(x_{1}\right),\ldots,\pi\left(x_{r}\right)\right)\right\Vert \leq C_{j}.\]
Choose $\delta_{c}$ and $\epsilon_{c}$ as before, with the $\delta_{c}$
summable. Keeping with our earlier notation, we are going to define
$\mathbf{z}^{(c)}$ from $\mathbf{z}^{(c-1)}$ by\[
z_{k}^{(c)}=\left(1-\delta_{c}u_{\lambda_{c}}\right)z_{k}^{(c-1)}.\]
Lemma~\ref{lem:approxRoots} and Theorem~\ref{thm:normTechnicalLemma}
give us 

\begin{align*}
 & \limsup_{\lambda}\left\Vert \sum_{k=1}^{r}\left(\left(1-\delta_{c}u_{\lambda}\right)z_{k}^{(c-1)}\left(z_{k}^{(c-1)}\right)^{*}\left(1-\delta_{c}u_{\lambda}\right)\right)^{\frac{q}{2}}\right\Vert \\
 & \hspace{1in}=\limsup_{\lambda}\left\Vert \sum_{k=1}^{r}\left(z_{k}^{(c-1)}\left(z_{k}^{(c-1)}\right)^{*}\right)^{\frac{q}{2}}\left(1-\delta_{c}u_{\lambda}\right)^{q}\right\Vert \\
 & \hspace{1in}\leq\limsup_{\lambda}\left\Vert \sum_{k=1}^{r}\left(z_{k}^{(c-1)}\left(z_{k}^{(c-1)}\right)^{*}\right)^{\frac{q}{2}}\left(1-\delta_{c}u_{\lambda}\right)\right\Vert \\
 & \hspace{1in}=1\end{align*}
and\begin{align*}
 & \limsup_{\lambda}\left\Vert p_{j}\left(\left(1-\delta_{1}u_{\lambda}\right)\mathbf{z}^{(c-1)}\right)\right\Vert \\
 & \hspace{1in}=\limsup_{\lambda}\left\Vert p_{j}\left(\mathbf{z}^{(c-1)}\right)\left(1-\delta_{c}u_{\lambda}\right)^{d_{j}}\right\Vert \\
 & \hspace{1in}\leq\limsup_{\lambda}\left\Vert p_{j}\left(\mathbf{z}^{(c-1)}\right)\left(1-\delta_{c}u_{\lambda}\right)\right\Vert \\
 & \hspace{1in}=C_{j}.\end{align*}
The limit of the $z_{j}^{(c)}$ will exist because \begin{align*}
\left\Vert z_{k}^{(c)}-z_{k}^{(c-1)}\right\Vert  & =\left\Vert \left(1-\delta_{c}u_{\lambda_{c}}\right)z_{k}^{(c-1)}-z_{k}^{(c-1)}\right\Vert \\
 & \leq\left\Vert \left(1-\delta_{c}u_{\lambda_{c}}\right)-1\right\Vert \left\Vert z_{k}^{(c-1)}\right\Vert \\
 & \leq\delta_{c}\left\Vert x_{k}\right\Vert .\end{align*}

\end{proof}

\section{Soft Cylinders\label{sec:Soft-Cylinders}}

When we stray from homogeneous relations, we come across $K$-theoretical
obstructions to projectivity. To illustrate what properties can still
hold, we offer the example of the {}``soft cylinder.'' The weaker
properties are semiprojectivity (as in \cite{Blackadar-shape-theory})
and the RFD property, meaning ``residually finite dimensional.'' Projectivity
implies semiprojectivity and also RFD (\cite[\S 1]{LorPederProjTrans}).

For $\epsilon\geq0$ we define the soft cylinder almost like Exel's
soft torus (\cite{ExelSoftTorusI}):\[
A_{\epsilon}=C_{1}^{*}\left\langle u,h\left|\,\strut u^{*}u=uu^{*}=1,\,-1\leq h\leq1,\,\left\Vert uh-hu\right\Vert \leq\epsilon\right.\right\rangle .\]
Notice we retained some homogeneity.
\begin{thm}
For positive $\epsilon,$ the soft cylinder $A_{\epsilon}$ is semiprojective.\end{thm}
\begin{proof}
Suppose $B$ is a unital $C^{*}$-algebra, with ideal $I=\overline{\bigcup I}_{n}$
for some increasing sequence of ideals $I_{n}.$ Suppose we are given
$u$ and $h$ in $B/I$ where $u$ is unitary, $-1\leq h\leq1$ and
\[
\left\Vert uh-hu\right\Vert \leq\epsilon.\]
For some $n$ it is possible to lift $u$ to $v$ in $B/I_{n}$ that
is a unitary (\cite[Prop.\ 2.21]{Blackadar-shape-theory}). Take any
lift of $h$ to $-1\leq k\leq1$ in $B/I_{n}.$ Theorem~\ref{thm:MixedVariables}
tells us there is $\hat{k}$ in $A/I$ with $-1\leq\hat{k}\leq1$
and\[
\left\Vert v\hat{k}-\hat{k}v\right\Vert \leq\epsilon.\]

\end{proof}
Eilers and Exel (\cite{EilersExelSoftTorusRFD}) have shown that the
soft torus is RFD. The same can be said, and proven much more easily,
for the soft cylinder.
\begin{thm}
For positive $\epsilon,$ the soft cylinder $A_{\epsilon}$ is RFD.\end{thm}
\begin{proof}
Consider the surjection\[
\rho:C(S^{1})\ast_{\mathbb{C}}C[0,1]\twoheadrightarrow\widetilde{A_{\epsilon}}\]
that sends the obvious unitary generator to $u$ and the obvious positive,
norm-one generator to $h.$ By \cite[Theorem 3.2]{ExelLoringRFDfreeProd}
the free product is RFD. Theorem~\ref{thm:MixedVariables} tells
us that $\rho$ is split. Thus $A_{\epsilon}$ can be embedded in
an RFD $C^{*}$-algebra and so is itself RFD.
\end{proof}
Our lifting theorems can be used to determine many more $C^{*}$-algebras
are RFD. The study of weak projectivity (\cite{LoringWeaklyProjective})
and RFD of the $C^{*}$-algebras associated to rather general relations
that have some homogeneity might lead to some interesting examples.
These topics will be explored elsewhere.

\section{Cones are Limits of Projective $C^{*}$-Algebras\label{sec:Cones-are-Limits}}

We end with a tantalizingly result: every cone is the limit of projectives.
As a $C^{*}$-algebra with a projective cone must be semiprojective
(\cite[II.8.3.10]{BlackadarOperatorAlgebras}) it would seem that
we are close to proving that every separable $C^{*}$-algebra is a
limit of semiprojective $C^{*}$-algebras. 

We say definitively, projectivity is not ``extremely rare'' (\cite[p.\ 73]{Loring-lifting-perturbing}).
\begin{lem}
\label{lem:relationsForCones} Suppose $A$ is the unital $C^{*}$-algebra
\[
A=C_{1}^{*}\left\langle x_{1},x_{2},\ldots\left|\begin{array}{c}
-C_{k}\leq x_{k}\leq C_{k}\quad(\forall k)\\
p_{j}\left(\mathbf{x}\right)=0\quad(\forall j)\end{array}\right.\right\rangle \]
where the $p_{1},p_{2}\ldots$ are NC polynomials in the $x_{k}$
of degrees $D_{j}$ with zero constant term. Then the cone $\mathbf{C}A$
has presentation \[
\mathbf{C}A=C^{*}\left\langle h,x_{1},x_{2},\ldots\left|\begin{array}{c}
0\leq h\leq1\\
hx_{k}=x_{k}h\quad(\forall k)\\
-C_{k}h\leq x_{k}\leq C_{k}h\quad(\forall k)\\
q_{j}\left(h,\mathbf{x}\right)=0\quad(\forall j)\end{array}\right.\right\rangle ,\]
where $q_{j}$ is the NC polynomial derived from the $p_{j}$ by padding
monomials on the left with various powers of $h$ so that $q_{j}$
is homogeneous with degree $D_{j}.$\end{lem}
\begin{proof}
To illustrate the construction of the $q_{j},$ if \[
p_{1}=x_{1}+3x_{1}x_{2}^{*}x_{1}\]
 then \[
q_{1}=h^{2}x_{1}+3x_{1}x_{2}^{*}x_{1}.\]
In general, we can break up $p_{j}$ into homogeneous summands \[
p_{j}=\sum_{d=1}^{D_{j}}p_{j,d}\]
and then describe the $q_{j}$ as \[
q_{j}\left(h,\mathbf{x}\right)=\sum_{d=1}^{D_{j}}h^{D_{j}-d}p_{j,d}\left(\mathbf{x}\right).\]

Let the universal $C^{*}$-algebra for these relations be denoted
$\mathcal{U}.$ This exists, as the relations satisfy the needed four
axioms as in \cite{LoringCstarRelations}. One of the axioms is that
setting all variables to the zero elements in $\{0\}$ leads to a
representation of the relations, which is true because we require
the constant terms to be zero. 

To define a $*$-homomorphism $\mathcal{U}\rightarrow\mathbf{C}A$
we define in $\mathbf{C}A=C_{0}\left((0,1],A\right)$ elements $\tilde{x}_{k}=tx_{k}$
and $\tilde{h}=t,$ shorthand for $\tilde{x}(t)=tx$ and so forth.
It is obvious that $0\leq\tilde{h}\leq1$ and that $\tilde{h}$ commutes
with each $\tilde{x}_{k}.$ Also \[
-C_{k}\leq x_{k}\leq C_{k}\implies-tC_{k}\leq tx_{k}\leq tC_{k}\]
so $-C_{j}\tilde{h}\leq\tilde{x}_{j}\leq C_{j}\tilde{h}.$ The last
relation holds as well since \begin{align*}
\left(q_{j}\left(\tilde{h},\tilde{\mathbf{x}}\right)\right)(t) & =\sum_{d=1}^{D_{j}}t^{D_{j}-d}p_{j,d}\left(t\mathbf{x}\right)\\
 & =t^{D_{j}}\sum_{d=1}^{D_{j}}p_{j,d}\left(\mathbf{x}\right)\\
 & =t^{D_{j}}p_{j}\left(\mathbf{x}\right)\\
 & =0.\end{align*}

Next we will show this map is onto. Basic algebra, and the usual isomorphism
of $\mathbf{C}A$ with $C_{0}(0,1]\otimes A,$ tells us that functions
of the form $t^{m}w$ generate the cone, where $w$ ranges over words
in the $x_{k}.$ Suppose $w=w_{1}w_{2}\cdots w_{n}.$ If $m\geq n$
then this is easily in the image, as \[
t^{m}w=t^{m-n}\left(tw_{1}\right)\left(tw_{2}\right)\ldots\left(tw_{n}\right).\]
If $1\leq m<n$ then the Stone-Weierstrass theorem tells us we can
approximate in $C_{0}(0,1]$ the function $t^{m}$ by a polynomial
in $t^{n},t^{n+1},\ldots$ and so can approximate $t^{m}w$ by a polynomial
in $t^{n}w,t^{n+1}w,\ldots$ and the map is indeed onto. We turn to
proving it is one-to-one. 

Consider an irreducible representation in $\mathbb{B}(\mathbb{H})$
of the relations defining $\mathcal{U}$ by $H$ and $X_{1},X_{2},\ldots.$
Since $0\leq H\leq1$ and $HX_{k}=X_{k}H$ and we find that $H$ is
central and so $H=\lambda I$ for some scalar $\lambda$ with $0\leq\lambda\leq1.$
If $\lambda=0$ then $H=0$ and \[
-C_{k}\lambda\leq X_{k}\leq C_{k}\lambda\implies X_{k}=0.\]
This is the zero representation, which is the pullback of the zero
representation of $\mathbf{C}A.$ If $\lambda$ is positive, then
\[
-C_{k}\lambda\leq X_{k}\leq C_{k}\lambda\implies-C_{k}\leq\lambda^{-1}X_{k}\leq C_{k}\]
and $q_{j}\left(H,\mathbf{X}\right)=0$ implies \begin{align*}
p_{j}\left(\lambda^{-1}\mathbf{X}\right) & =\sum_{d=1}^{D_{j}}\lambda^{-d}p_{j,d}\left(\mathbf{X}\right)\\
 & =\lambda^{-D_{j}}\sum_{d=1}^{D_{j}}H^{D_{j}-d}p_{j,d}\left(\mathbf{X}\right)\\
 & =\lambda^{-D_{j}}q_{j}\left(H,\mathbf{X}\right)\\
 & =0.\end{align*}
Thus the $\lambda^{-1}X_{k}$ form a representation of $A$ on $\mathbb{H}$
and so a representation of $\mathbf{C}A$ via the composition \input{insert3.tex}
This sends $\tilde{h}$ to $\lambda I=H$ and $\tilde{x}_{k}$ to
$X_{k},$ finishing the proof.\end{proof}
\begin{thm}
\label{thm:ConeLikeProjectives} If $q_{1},q_{2}\ldots$ are homogeneous
NC polynomials, each of degree at least one, in noncommuting variables
$h,x_{1},x_{2},\ldots$ then for positive constants $C_{1},\ldots,C_{J}$
and $D_{1},\ldots,D_{K},$ the $C^{*}$-algebra \[
C^{*}\left\langle h,x_{1},x_{2},\ldots\left|\begin{array}{c}
0\leq h\leq1,\\
-C_{k}h\leq x_{k}\leq C_{k}h,\quad(1\leq k\leq K)\\
\left\Vert q_{j}\left(h,x_{1},x_{2},\ldots\right)\right\Vert \leq D_{j}\quad(1\leq j\leq J)\end{array}\right.\right\rangle \]
is projective.\end{thm}
\begin{proof}
We can find some $r$ so that $x_{r+1},x_{r+2},\ldots$ are not in
any of the polynomials $q_{1},\ldots,q_{J}.$ If we relabel these
$y_{1},y_{2},\ldots$ our lifting problem becomes\[
0\leq h\leq1,\]
\[
-C_{k}h\leq x_{k}\leq C_{k}h,\]
\[
-C_{k}^{\prime}h\leq y_{k}\leq C_{k}^{\prime}h,\]
\[
\left\Vert q_{j}\left(h,\mathbf{x}\right)\right\Vert \leq D_{j}\quad(1\leq j\leq J)\]
 where now the $q_{j}$ are homogeneous in $\left\{ h,x_{1},\ldots,x_{r}\right\} .$
We are using $\mathbf{x}$ for $\left(x_{1},\ldots,x_{r}\right).$ 

Given $h,$ $ $$x_{k}$ and $y_{k}$ in $A$ with\[
0\leq\pi\left(h\right)\leq1,\]
\[
-C_{k}\pi\left(h\right)\leq\pi\left(x_{k}\right)\leq C_{k}\pi\left(h\right),\]
\[
-C_{k}^{\prime}\pi\left(h\right)\leq\pi\left(y_{k}\right)\leq C_{k}^{\prime}\pi\left(h\right),\]
\[
\left\Vert q_{j}\left(\pi\left(h\right),\pi\left(\mathbf{x}\right)\right)\right\Vert \leq D_{j}\quad(1\leq j\leq J)\]
we first find a new lift $\hat{h}$ of $\pi(h)$ with \[
0\leq\hat{h}\leq1.\]
Using Davidson's two-sided order lifting theorem (\cite{Davidson-Lifting-positive})
we find $\hat{x}_{k}$ and $\hat{y}_{k}$ with\[
-C_{k}\hat{h}\leq\hat{x}_{k}\leq C_{k}\hat{h},\]
\[
-C_{k}^{\prime}\hat{h}\leq\hat{y}_{k}\leq C_{k}^{\prime}\hat{h}\]
and $\pi\left(\hat{x}_{k}\right)=\pi\left(x_{j}\right)$ and $\pi\left(\hat{y}_{k}\right)=\pi\left(y_{j}\right).$
By Theorem \ref{thm:MixedVariables} there is an $m$ with $0\leq m\leq1$
in $1+I$ so that\[
\left\Vert q_{j}\left(m\hat{h}m,m\hat{\mathbf{x}}m\right)\right\Vert \leq D_{j}.\]
Our desired lifts are $ $$m\hat{h}m,$ $m\hat{\mathbf{x}}m$ and
$m\hat{\mathbf{y}}m.$ \end{proof}
\begin{lem}
\label{lem:C*fromPolyRelations} Let $D$ be a separable $C^{*}$-algebra.
Then\[
D\cong C^{*}\left\langle x_{1},x_{2},\ldots\left|\begin{array}{c}
-C_{j}\leq x_{j}\leq C_{j}\quad(\forall j)\\
p_{k}\left(x_{1},x_{2},\ldots\right)=0\quad(\forall k)\end{array}\right.\right\rangle \]
for a countable collection of NC polynomials. \end{lem}
\begin{proof}
Example 1.3(b) in \cite{Blackadar-shape-theory} tells us that $D$
has a presentation with countably many generators, countably many
relations in the form of a NC $*$-polynomial set to zero and countably
many norm conditions. We will modify Blackadar's method a bit.

Let $\mathbb{F}=\mathbb{Q}+i\mathbb{Q},$ which is a countable dense
subfield of $\mathbb{C}.$ Select a countable dense sequence in $D$
and apply to this sequence all polynomials over $\mathbb{F}$ in countably
many variables. This results in a countable, dense $\mathbb{F}$-$*$-subalgebra
$B$ of $D.$ Enumerate $B$ as $x_{1},x_{2},\ldots.$ The algebraic
operations for $B$ can be encoded in $*$-polynomial relations. For
example, if $\alpha x_{j}=x_{k}$ for some $\alpha$ in $\mathbb{F},$
then we use the relation $\alpha x_{j}-x_{k}=0.$ If $x_{j}^{*}=x_{k}$
then we use the relation $x_{j}^{*}-x_{k}=0,$ and so forth. This
means $B$ is the universal $\mathbb{F}$-$*$-algebra for generators
$x_{1},x_{2},\ldots$ and some countable set of $*$-polynomial relations
$p_{j}(x_{1},x_{2},\ldots)=0.$ We now add to these relations the
$C^{*}$-relations $\left\Vert x_{k}\right\Vert \leq C_{k}$ where
$C_{k}$ is the norm of the element $x_{k}$ in $D.$ Then any function
$f:B\rightarrow G,$ for $G$ a $C^{*}$-algebra, that satisfies these
relations is first of all an $\mathbb{F}$-linear $*$-algebra homomorphism.
It is continuous with respect to the norm on $D$ since $x_{j}-x_{k}$
will equal some $x_{\ell}$ so we have the relation $\left\Vert f(x_{\ell})\right\Vert \leq\left\Vert x_{\ell}\right\Vert $
and so \[
\left\Vert f(x_{j})-f(x_{k})\right\Vert =\left\Vert f(x_{\ell})\right\Vert \leq\left\Vert x_{\ell}\right\Vert =\left\Vert x_{j}-x_{k}\right\Vert .\]
It therefore extends to a continuous function $\varphi:D\rightarrow G.$
This extended function will be linear over $\mathbb{C}.$ To verify
this, consider $\alpha=\lim\alpha_{n},$ a limit of scalars from $\mathbb{F},$
and $d=\lim d_{j},$ a limit of elements in $B.$ Then \[
\varphi(\alpha_{n}d)=\varphi(\lim_{j}\alpha_{n}d_{j})=\lim_{j}f(\alpha_{n}d_{j})=\alpha_{n}\lim_{j}f(d_{j})=\alpha_{n}\varphi(d)\]
and \[
\varphi(\alpha d)=\varphi(\lim_{n}\alpha_{n}d)=\lim_{n}\varphi(\alpha_{n}d)=\lim_{n}\alpha_{n}\varphi(d)=\alpha\varphi(d).\]
Finally, continuity implies that $\varphi$ Is a $*$-homomorphism.
It is uniquely determined by $f$ and so $D$ is universal for these
relations.

We can eliminate many of the norm conditions. Suppose we keep only
the norm restrictions $\left\Vert x_{k}\right\Vert \leq C_{k}$ for
those $x_{k}$ that are self-adjoint. Then the estimate that gave
continuity changes a little. Any $x_{j}-x_{k}$ will equal some $x_{\ell}$
and for some $r$ and $s$ we will have $x_{r}=\frac{1}{2}x_{\ell}+\frac{1}{2}x_{\ell}^{*}$
and $x_{s}=\frac{-i}{2}x_{\ell}-\frac{i}{2}x_{\ell}^{*}.$ Of course
these are the real and imaginary part of $x_{\ell},$ and as they
are self-adjoint we have the relations $\left\Vert f(x_{r})\right\Vert \leq\left\Vert x_{r}\right\Vert $
and $\left\Vert f(x_{s})\right\Vert \leq\left\Vert x_{s}\right\Vert .$
Therefore \begin{align*}
\left\Vert f(x_{j})-f(x_{k})\right\Vert  & =\left\Vert f(x_{r})+if(x_{s})\right\Vert \\
 & \leq\left\Vert f(x_{r})\right\Vert +\left\Vert f(x_{s})\right\Vert \\
 & \leq\left\Vert x_{r}\right\Vert +\left\Vert x_{s}\right\Vert \\
 & \leq2\left\Vert x_{\ell}\right\Vert \\
 & =2\left\Vert x_{j}-x_{k}\right\Vert .\end{align*}
This still gives us continuity and so the rest of the proof goes through.

We can toss the generators that are not self-adjoint if we modify
each polynomial by the evaluating $x_{k}$ at $x_{r}+ix_{s}$ whenever
$x_{r}$ and $x_{s}$ are the real and imaginary parts of $x_{k}.$
Among the polynomial relations will be $x_{j}^{*}-x_{j}=0$ for the
generators we are keeping. Given this, it is our option to use the
relation $\left\Vert x_{j}\right\Vert \leq C_{j}$ or $-C_{j}\leq x_{j}\leq C_{j}.$\end{proof}
\begin{thm}
If $A$ is a separable $C^{*}$-algebra then its cone $\mathbf{C}A$
is isomorphic to the inductive limit of a countable system of projective
$C^{*}$-algebras with surjective bonding maps. \end{thm}
\begin{proof}
We start with the case where $A=\widetilde{D}$ for some separable,
possibly unital $C^{*}$-algebra. 

Lemma~\ref{lem:C*fromPolyRelations} tells us \[
A\cong C_{1}^{*}\left\langle x_{1},x_{2},\ldots\left|\begin{array}{c}
-C_{k}\leq x_{k}\leq C_{k}\quad(\forall k)\\
p_{j}\left(x_{1},x_{2},\ldots\right)=0\quad(\forall j)\end{array}\right.\right\rangle \]
and then Lemma~\ref{lem:relationsForCones} tells us \[
\mathbf{C}A\cong C^{*}\left\langle h,x_{1},x_{2},\ldots\left|\begin{array}{c}
0\leq h\leq1\\
hx_{k}=x_{k}h\quad(\forall k)\\
-C_{k}h\leq x_{k}\leq C_{k}h\quad(\forall k)\\
q_{j}\left(h,x_{1},x_{2},\ldots\right)=0\quad(\forall j)\end{array}\right.\right\rangle \]
where the $q_{k}$ are homogeneous. Clearly \[
\mathbf{C}A\cong\lim_{\longrightarrow}P_{n}\]
where \[
P_{n}=C^{*}\left\langle h,x_{1},x_{2},\ldots\left|\begin{array}{c}
0\leq h\leq1\\
-C_{k}h\leq x_{k}\leq C_{k}h,\quad(\forall k)\\
\left\Vert hx_{k}-x_{k}h\right\Vert \leq\frac{1}{n},\quad(k=1,\ldots,n)\\
\left\Vert q_{j}\left(x_{1},x_{2},\ldots\right)\right\Vert \leq\frac{1}{n}\quad(j=1,\ldots,n)\end{array}\right.\right\rangle .\]
Since the commutators are homogeneous NC polynomials, Theorem~\ref{thm:ConeLikeProjectives}
applies and the $P_{n}$ are projective. We are done for $\mathbf{C}A=\mathbf{C}\left(\widetilde{D}\right).$
What about $\mathbf{C}D$?

We have the exact sequence \input{insert4.tex} Of course $\mathbf{C\mathbb{C}}$
equals $C_{0}(0,1]$ and is projective. Let $Q_{n}$ be the kernel
of the map of $P_{n}$ onto $C_{0}(0,1]$ that sends $h$ to $t\mapsto t$
and $x_{k}$ to zero. Then we have \input{insert5.tex} with the rows
exact. Also, $\mathbf{C}D$ is isomorphic to ${\displaystyle \lim_{\longrightarrow}Q_{n}},$
which we can see as follows. 

There is a $*$-homomorphism \[
\varphi:\lim_{\longrightarrow}Q_{n}\rightarrow\mathbf{C}D\]
 induced by the maps $Q_{n}\rightarrow\mathbf{C}D.$ The maps $Q_{n}\rightarrow P_{n}$
are inclusions and hence isometries. Theorem~13.1.2.2 in \cite{Loring-lifting-perturbing}
implies that the induced map ${\displaystyle \lim_{\longrightarrow}Q_{n}\rightarrow\lim_{\longrightarrow}P_{n}}$
is an also an isometry. From the commutative diagram \input{insert6.tex}
we conclude $\varphi$ is injective. As to surjectivity, consider
$x$ in $\mathbf{C}D.$ This gets sent to $0$ in $\mathbf{C}\mathbb{C}.$
Any lift of $x$ to $y$ in $P_{1}$ is also sent to zero in $\mathbf{C}\mathbb{C},$
and so is $Q_{n}.$ This shows $\varphi$ is onto.

By Theorem~5.3 of \cite{LorPederProjTrans} the $Q_{n}$ are projective. 
\end{proof}

\end{document}

%% file: insert3.tex
\[
\xymatrix{
\mathbf C A \ar[r]^{\delta_\lambda}& A \ar[r] & \mathbb{B}(\mathbb{H}) .
}
\]

%% file: insert4.tex
\[
\xymatrix{
0 \ar[r] 
	& \mathbf C D \ar[r] 
		&  \mathbf C A \ar[r]  
			& \mathbf C \mathbb C \ar[r] 
				& 0.
}
\]

%% file: insert5.tex
\[
\xymatrix{
0 \ar[r] 
	& Q_n \ar[r] \ar[d] 
		&  P_n \ar[r] \ar[d]  
			& \mathbf C \mathbb C \ar[r] \ar@{=}[d] 
				& 0
\\
0 \ar[r] 
	& Q_{n+1} \ar[r] \ar[d] 
		&  P_{n+1} \ar[r] \ar[d]  
			& \mathbf C \mathbb C \ar[r] \ar@{=}[d] 
				& 0
\\
0 \ar[r] 
	& \mathbf C D \ar[r] 
		&  \mathbf C A \ar[r]  
			& \mathbf C \mathbb C \ar[r] 
				& 0
}
\]

%% file: insert6.tex
\[
\xymatrix{
\displaystyle{\lim_{\longrightarrow} Q_n } 
				\ar@{^{(}->}[r]  \ar[d] ^{\varphi} 
	& \displaystyle{\lim_{\longrightarrow} P_n } \ar[d]^{\cong}  
\\
\mathbf C D \ar@{^{(}->}[r] 
		&  \mathbf C A 
}
\]